\documentclass[twoside]{article}
\usepackage{mathrsfs}
\usepackage{amsmath}
\usepackage{amsfonts}
\usepackage{amssymb}
\usepackage{amsmath,amssymb,amsthm}
\usepackage{amsmath,amssymb,amsthm,amscd}
\usepackage[frame,cmtip,arrow,matrix,line,graph,curve]{xy}
\usepackage{graphpap, color}
\usepackage[mathscr]{eucal}
\usepackage{color}
\usepackage{verbatim}
\usepackage{hyperref}
\usepackage{cite}
\numberwithin{equation}{section} \numberwithin{equation}{section}

\renewcommand{\oddsidemargin}{5mm}

\theoremstyle{plain}\newtheorem{thm}{Theorem}[section]
\theoremstyle{plain}\newtheorem{lem}{Lemma}[section]
\theoremstyle{definition}
\theoremstyle{definition}
\theoremstyle{plain}
\theoremstyle{plain}
\theoremstyle{remark}\newtheorem{rem}{Remark}[section]

\usepackage{fancyhdr}
\begin{document}

\title{Estimates for Eigenvalues of Poly-harmonic Operators
}
\author{Guoxin Wei~~and~~Lingzhong Zeng}

\vskip 5mm

\medskip
\date{}
\maketitle

\abstract {In this paper, we study eigenvalues of the poly-Laplacian
with arbitrary order on a bounded domain in an $n$-dimensional
Euclidean space and obtain a lower bound for eigenvalues, which
generalizes the results due to Cheng-Wei \cite{CW2} and gives an
improvement of results due to Cheng-Qi-Wei
\cite{CQW}.}\\



\normalsize
\section{Introduction}\label{sec1}

Let $\Omega$ be a bounded domain with
piecewise smooth boundary $\partial\Omega$ in an $n$-dimensional
Euclidean space $\mathbb{R}^{n}$. Let $\lambda_{i}$ be the $i$-th
eigenvalue of Dirichlet eigenvalue problem of the poly-Laplacian
with arbitrary order:
\begin{equation}\label{1.1}\begin{cases}(-\Delta)^{l}u = \lambda u,~~~~~~~~~~~~~~~~~~~~~~~~~~~~~~~~~~~~~~in~\Omega,\\[2mm]
u=\dfrac{\partial u}{\partial\nu}=\cdots=\dfrac{\partial^{l-1}
                  u}{\partial\nu^{l-1}}=0,~~~~~~~~~~~~~~~~~~~~on~\partial\Omega,\end{cases}\end{equation}
where $\Delta$ is the Laplacian in $\mathbb{R}^{n}$ and $\nu$
denotes the outward unit normal vector field of the boundary
$\partial\Omega$. It is well known that the spectrum of this
eigenvalue problem (\ref{1.1}) is real and discrete:
$$0<\lambda_{1}\leq\lambda_{2}\leq\lambda_{3}\leq\cdots\rightarrow+\infty,$$ where each $\lambda_{i}$ has finite
multiplicity which is repeated according to its multiplicity. Let
$V(\Omega)$ denote the volume of $\Omega$ and let $B_{n}$ denote the
volume of the unit ball in $\mathbb{R}^{n}$.

When $l = 1$, the eigenvalue problem (\ref{1.1}) is called a fixed
membrane problem. In this case, one has the following Weyl's
asymptotic formula
\begin{equation}\label{1.2}\lambda_{k}\sim \frac{4\pi^{2}}{(B_{n}V(\Omega))^{\frac{2}{n}}}k^{\frac{2}{n}},~~k \rightarrow+\infty.\end{equation}
From the above asymptotic formula, one can derive
\begin{equation}\label{1.3}\frac{1}{k}\sum^{k}_{i=1}\lambda_{i}\sim\frac{n}{n+ 2}\frac{4\pi^{2}}
{(B_{n}V(\Omega))^{\frac{2}{n}}}k^{\frac{2}{n}},~~k
\rightarrow+\infty.
\end{equation} P$\acute{\textnormal{o}}$lya \cite{Pol} proved that \begin{equation}\label{1.4}\lambda_{k}\geq \frac{4\pi^{2}}
{(B_{n}V(\Omega))^{\frac{2}{n}}}k^{\frac{2}{n}},~~{\rm for}\ k
=1,2,\cdots,\end{equation} if $\Omega$ is a tiling domain in
$\mathbb{R}^{n}$ . Furthermore, he proposed a conjecture as follows:\\

\noindent\textbf{Conjecture of P$\acute{\textnormal{o}}$lya}. \emph{If
$\Omega$ is a bounded domain in $\mathbb{R}^{n}$, then the $k$-th
eigenvalue $\lambda_{k}$ of the fixed membrane problem satisfies}
\begin{equation}\label{1.5}\lambda_{k}\geq \frac{4\pi^{2}}
{(B_{n}V(\Omega))^{\frac{2}{n}}}k^{\frac{2}{n}},~~{\rm for}\ k
=1,2,\cdots.\end{equation} On the conjecture of
P$\acute{\textnormal{o}}$lya, Berezin \cite{Be} and Lieb \cite{Lie}
gave a partial solution.  In particular, Li and Yau \cite{LiY}
proved that
\begin{equation}\label{1.6}\frac{1}{k}\sum^{k}_{i=1}\lambda_{i}\geq\frac{n}{n+
2}\frac{4\pi^{2}}
{(B_{n}V(\Omega))^{\frac{2}{n}}}k^{\frac{2}{n}},~~\textnormal{for}~
k =1,2,\cdots.
\end{equation} The formula (\ref{1.3}) shows that the result
of Li and Yau is sharp in the sense of average. From this formula
(\ref{1.6}), one can infer
\begin{equation}\label{1.7}\lambda_{k}\geq\frac{n}{n+
2}\frac{4\pi^{2}}
{(B_{n}V(\Omega))^{\frac{2}{n}}}k^{\frac{2}{n}},~~\textnormal{for}~
k =1,2,\cdots,
\end{equation} which gives a partial solution for the
conjecture of P$\acute{\textnormal{o}}$lya with a factor
$\dfrac{n}{n+2}$. Recently, Melas \cite{M} has improved the estimate
(\ref{1.6}) to the following:
\begin{equation}\label{1.8}\frac{1}{k}\sum^{k}_{i=1}\lambda_{i}\geq\frac{n}{n+
2}\frac{4\pi^{2}}
{(B_{n}V(\Omega))^{\frac{2}{n}}}k^{\frac{2}{n}}+\frac{1}{24(n+2)}\frac{V(\Omega)}{I(\Omega)},~~\textnormal{for}~k
=1,2,\cdots,
\end{equation}
where
$$I(\Omega)=\min_{a\in\mathbb{R}^{n}}\int_{\Omega}|x-a|^{2}dx$$ is
called \emph{the moment of inertia} of $\Omega$.

When $l = 2$, the eigenvalue problem (\ref{1.1}) is called a clamped
plate problem. For the eigenvalues of the clamped plate problem,
Agmon \cite{A} and Pleijel \cite{Pl} obtained
\begin{equation}\label{1.9}\lambda_{k}\sim
\frac{16\pi^{4}}{(B_{n}V(\Omega))^{\frac{4}{n}}}k^{\frac{4}{n}},~~k
\rightarrow+\infty.\end{equation} From the above formula
(\ref{1.9}), one can obtain
\begin{equation}\label{1.10}\frac{1}{k}\sum^{k}_{i=1}\lambda_{i}\sim\frac{n}{n+4}\frac{16\pi^{4}}{(B_{n}V(\Omega))^{\frac{4}{n}}}k^{\frac{4}{n}},~~k
\rightarrow+\infty.
\end{equation}Furthermore, Levine and Protter \cite{LP}
proved that the eigenvalues of the clamped plate problem satisfy the
following inequality:
\begin{equation}\label{1.11}\frac{1}{k}\sum^{k}_{i=1}\lambda_{i}\geq\frac{n}{n+4}\frac{16\pi^{4}}{(B_{n}V(\Omega))^{\frac{4}{n}}}k^{\frac{4}{n}}.
\end{equation} The formula (\ref{1.10}) shows that the coefficient of $k^{\frac{4}{n}}$ is the best possible constant.
By adding to its right hand side two terms of lower order in $k$,
Cheng and Wei \cite{CW1} obtained the following estimate which is an
improvement of (\ref{1.11}):
\begin{eqnarray}\label{1.12}\begin{aligned}\frac{1}{k}\sum^{k}_{i=1}\lambda_{i}&\geq\frac{n}{n+
4}\frac{16\pi^{4}} {(B_{n}V(\Omega))^{\frac{4}{n}}}k^{\frac{4}{n}}
\\&~\quad+\Bigg{(}\frac{n+2}{12n(n+4)}-\frac{1}{1152n^{2}(n+4)}\Bigg{)}\frac{4\pi^{2}}{(B_{n}V(\Omega))^{\frac{2}{n}}}\frac{n}{n+2}\frac{V(\Omega)}{I(\Omega)}k^{\frac{2}{n}}
\\&~\quad+\Bigg{(}\frac{1}{576n(n+4)}-\frac{1}{27648n^{2}(n+2)(n+4)}\Bigg{)}\Bigg{(}\frac{V(\Omega)}{I(\Omega)}\Bigg{)}^{2}.
\end{aligned}\end{eqnarray}
Very recently, Cheng and Wei \cite{CW2} have improved the estimate
(\ref{1.12}) to the following:
\begin{eqnarray}\label{1.13}\begin{aligned}\frac{1}{k}\sum^{k}_{i=1}\lambda_{i}&\geq\frac{n}{n+
4}\frac{16\pi^{4}}
{(B_{n}V(\Omega))^{\frac{4}{n}}}k^{\frac{4}{n}}\\&~\quad+\frac{n+2}{12n(n+4)}\frac{4\pi^{2}}
{(B_{n}V(\Omega))^{\frac{2}{n}}}\frac{n}{n+
2}\frac{V(\Omega)}{I(\Omega)}k^{\frac{2}{n}}\\&~\quad+\frac{(n+2)^{2}}{1152n(n+4)^{2}}\Bigg{(}\frac{V(\Omega)}{I(\Omega)}\Bigg{)}^{2}.
\end{aligned}\end{eqnarray}

When $l$ is arbitrary, Levine and Protter \cite{LP} proved the
following
\begin{equation}\label{1.14}\frac{1}{k}\sum^{k}_{i=1}\lambda_{i}\geq\frac{n}{n+2l}\frac{\pi^{2l}}{(B_{n}V(\Omega))^{\frac{2l}{n}}}k^{\frac{2l}{n}}, ~~\textnormal{for}~~k
=1,2,\cdots,
\end{equation}
which implies that
\begin{equation}\label{1.15}\lambda_{k}\geq\frac{n}{n+2l}\frac{\pi^{2l}}{(B_{n}V(\Omega))^{\frac{2l}{n}}}k^{\frac{2l}{n}},
~~\textnormal{for}~~k =1,2,\cdots.
\end{equation}
By adding $l$ terms of lower order of $k^{\frac{2l}{n}}$ to its
right hand side, Cheng, Qi and Wei \cite{CQW} obtained more sharper
result than (\ref{1.14}):

\begin{eqnarray}\label{1.16}\begin{aligned}\frac{1}{k}\sum^{k}_{i=1}\lambda_{i}&\geq\frac{n}{n+
2l}\frac{(2\pi)^{2l}}
{(B_{n}V(\Omega))^{\frac{2l}{n}}}k^{\frac{2l}{n}}+\frac{n}{(n+2l)}\\&~\quad\times\sum_{p=1}^{l}
\frac{l+1-p}{(24)^{p}n\cdots(n+2p-2)}\frac{(2\pi)^{2(l-p)}}{(B_{n}V(\Omega))^{\frac{2(l-p)}{n}}}
\Bigg{(}\frac{V(\Omega)}{I(\Omega)}\Bigg{)}^{p}k^{\frac{2(l-p)}{n}}.
\end{aligned}\end{eqnarray}

In this paper, we investigate eigenvalues of the Dirichlet
eigenvalue problem (\ref{1.1}) of Laplacian with arbitrary order and
prove the following:

\begin{thm}\label{thm1.1} Let $\Omega$ be a bounded domain in an n-dimensional
Euclidean space $\mathbb{R}^{n}$. Assume that $l\geq2$ and
$\lambda_{i}$ is the $i$-th eigenvalue of the eigenvalue problem
\textnormal{(\ref{1.1})}. Then the eigenvalues satisfy
\begin{eqnarray}\label{1.17}\begin{aligned}\frac{1}{k}\sum_{j=1}^{k}\lambda_{j}&\geq
\dfrac{n}{n+ 2l}\dfrac{(2\pi)^{2l}}
{(B_{n}V(\Omega))^{\frac{2l}{n}}}k^{\frac{2l}{n}}\\
&~\quad+\dfrac{l}{24(n+2l)}\dfrac{(2\pi)^{2(l-1)}}{(B_{n}V(\Omega))^{\frac{2(l-1)}{n}}}
\dfrac{V(\Omega)}{I(\Omega)}k^{\frac{2(l-1)}{n}}\\
&~\quad+\dfrac{l(n+2(l-1))^{2}}{2304n(n+2l)^{2}}
\dfrac{(2\pi)^{2(l-2)}}{(B_{n}V(\Omega))^{\frac{2(l-2)}{n}}}
\Bigg{(}\frac{V(\Omega)}{I(\Omega)}\Bigg{)}^{2}k^{\frac{2(l-2)}{n}}.
\end{aligned}\end{eqnarray}
\end{thm}

\begin{rem}\label{rem1.1} \emph{When $l =2$, Theorem
{\rm \ref{thm1.1}} reduces to the result of Cheng-Wei {\rm
\cite{CW2}}.}
\end{rem}

\begin{rem}\label{rem1.2} \emph{When $l\geq2$, we give an important improvement of the result {\rm (\ref{1.16})} due
to Cheng-Qi-Wei {\rm \cite{CQW}} since the inequality {\rm
(\ref{1.17})} is sharper than the inequality {\rm (\ref{1.16})}.
About this fact, we will give a proof in Section {\rm \ref{sec3}}.}
\end{rem}

\section{A Key Lemma}\label{sec2}

In this section, we will give a key Lemma which will play an
important role in the proof of Theorem \ref{thm1.1}.

\begin{lem}\label{lem2.1}
Let $b\geq2$ be a positive real number and $\mu>0$. If $\psi:
[0,~+\infty)\rightarrow[0,~+\infty)$ is a decreasing function such
that $$-\mu\leq\psi^{\prime}(s)\leq0$$ and
$$A:=\int_{0}^{\infty}s^{b-1}\psi(s)ds>0,$$ then, for any positive integer $l\geq2$, we have
\begin{eqnarray}\label{2.1}\begin{aligned}\int^{\infty}_{0}s^{b+2l-1}\psi(s)ds&\geq
\frac{1}{b+2l}(bA)^{\frac{b+2l}{b}}\psi(0)^{-\frac{2l}{b}}\\&~\quad+\frac{l}{6b(b+2l)\mu^{2}}(bA)^{\frac{b+2(l-1)}{b}}\psi(0)^{\frac{2b-2l+2}{b}}
\\&~\quad+\frac{l(b+2(l-1))^{2}}{144b^{2}(b+2l)^{2}\mu^{4}}(bA)^{\frac{b+2l-4}{b}}\psi(0)^{\frac{4b-2l+4}{b}}.\end{aligned}\end{eqnarray}
\end{lem}

\noindent\emph{Proof.} Let
\begin{eqnarray}\label{2.2}\varrho(t)=\frac{\psi\big{(}\frac{\psi(0)}{\mu}t\big{)}}{\psi(0)},\end{eqnarray}
then we have $\varrho(0)=1$ and $-1\leq\varrho'(t)\leq0.$ Without
loss of generality, we can assume
$$\psi(0)=1~\textnormal{and}~\mu=1.$$ Define
$$D_{l}:=\int^{\infty}_{0}s^{b+2l-1}\psi(s)ds.$$
One can assume that $D_{l}<\infty$, otherwise there is nothing to prove. Since $D_{l}<\infty$,
 we can conclude that
$$\lim_{s\rightarrow\infty}s^{b+2l-1}\psi(s)=0.$$ Putting
$h(s)=-\psi'(s)$ for $s\geq0$, we get $$0\leq
h(s)\leq1~~\textnormal{and}~~\int^{\infty}_{0}h(s)ds=\psi(0)=1.$$ By
making use of integration by parts, one has
\begin{eqnarray*}\int^{\infty}_{0}s^{b}h(s)ds=b\int^{\infty}_{0}s^{b-1}\psi(s)ds=bA,\end{eqnarray*}and
\begin{eqnarray*}\int^{\infty}_{0}s^{b+2l}h(s)ds\leq(b+2l)D_{l},\end{eqnarray*} since
$\psi(s)>0.$ By the same assertion as in \cite{M}, one can infer
that there exists an $\epsilon\geq0$ such that
\begin{eqnarray}\label{2.3}\int^{\epsilon+1}_{\epsilon}s^{b}ds=\int^{\infty}_{0}s^{b}h(s)ds=bA,\end{eqnarray}
and
\begin{eqnarray}\label{2.4}\int^{\epsilon+1}_{\epsilon}s^{b+2l}ds\leq\int^{\infty}_{0}s^{b+2l}h(s)ds\leq(b+2l)D_{l}.\end{eqnarray}
Let
$$\Theta(s)=bs^{b+2l}-(b+2l)\tau^{2l}s^{b}+2l\tau^{b+2l}-2l\tau^{b+2(l-1)}(s-\tau)^{2},$$
then we can prove that $\Theta(s)\geq0.$ By integrating the function
$\Theta(s)$ from $\epsilon$ to $\epsilon+1$, we deduce from
(\ref{2.3}) and (\ref{2.4}), for any $\tau>0,$
\begin{eqnarray}\label{2.5}b(b+2l)D_{l}-(b+2l)\tau^{2l}bA+2l\tau^{b+2l}\geq\frac{l}{6}\tau^{b+2(l-1)}.\end{eqnarray}
Define
$$f(\tau):=(b+2l)\tau^{2l}bA-2l\tau^{b+2l}+\frac{l}{6}\tau^{b+2(l-1)},$$
then we can obtain from (\ref{2.5}) that, for any $\tau>0,$
$$D_{l}=\int^{\infty}_{0}s^{b+2l-1}\psi(s)ds\geq\frac{f(\tau)}{b(b+2l)}.$$
Taking
$$\tau=(bA)^{\frac{1}{b}}\Bigg{(}1+\frac{b+2(l-1)}{12(b+2l)}(bA)^{-\frac{2}{b}}\Bigg{)}^{\frac{1}{b}},$$
then one has
\begin{eqnarray}\label{2.6}\begin{aligned}f(\tau)&=(bA)^{\frac{b+2l}{b}}
\Bigg{(}b-\frac{l(b+2(l-1))}{6(b+2l)}(bA)^{-\frac{2}{b}}\Bigg{)}\Bigg{(}1+\frac{b+2(l-1)}{12(b+2l)}(bA)^{-\frac{2}{b}}
\Bigg{)}^{\frac{2l}{b}}\\&~\quad+
\frac{l}{6}(bA)^{\frac{b+2(l-1)}{b}}\Bigg{(}1+\frac{b+2(l-1)}{12(b+2l)}(bA)^{-\frac{2}{b}}\Bigg{)}^{\frac{b+2(l-1)}{b}}.\end{aligned}\end{eqnarray}

Next, we consider four cases:

\noindent\textbf{Case 1: $\textbf{\emph{b}}\geq\textbf{2\emph{l}.}$}
For $t>0$, we have from the Taylor formula
\begin{eqnarray*}\begin{aligned}(1+t)^{\frac{2l}{b}}&\geq1+\frac{2l}{b}t+\frac{2l(2l-b)}{2b^{2}}t^{2}+
\frac{2l(2l-b)(2l-2b)}{6b^{3}}t^{3}\\&~\quad+\frac{2l(2l-b)(2l-2b)(2l-3b)}{24b^{4}}t^{4}
\end{aligned}\end{eqnarray*}
and
\begin{eqnarray*}\begin{aligned}(1+t)^{\frac{b+2(l-1)}{b}}&\geq1+\frac{2(l-1)+b}{b}t+\frac{(2(l-1)+b)(l-1)}{b^{2}}t^{2}
\\&~\quad+\frac{(2(l-1)-b)(l-1)(2(l-1)+b)}{3b^{3}}t^{3}.\end{aligned}\end{eqnarray*}
Putting $t=\frac{b+2(l-1)}{12(b+2l)}(bA)^{-\frac{2}{b}}$, we have
from
$(bA)^{\frac{2}{b}}\geq\frac{1}{(b+1)^{\frac{2}{b}}}\geq\frac{1}{3}>\frac{1}{4}$
(also see \cite{CW2}) that $t<\frac{1}{3}$ and
$b-2lt>\frac{4l}{3}>0$. And then, we obtain
\begin{eqnarray*}\begin{aligned}&~\quad\Bigg{(}b-\frac{l(b+2(l-1))}{6(b+2l)}(bA)^{-\frac{2}{b}}\Bigg{)}
\Bigg{(}1+\frac{b+2(l-1)}{12(b+2l)}(bA)^{-\frac{2}{b}}\Bigg{)}^{\frac{2l}{b}}\\&=
(b-2lt)(1+t)^{\frac{2l}{b}}\end{aligned}\end{eqnarray*}
\begin{eqnarray*}\begin{aligned}&\geq(b-2lt)\Bigg{[}1+\frac{2l}{b}t+\frac{2l(2l-b)}{2b^{2}}t^{2}+
\frac{2l(2l-b)(2l-2b)}{6b^{3}}t^{3}\\&~\quad+\frac{2l(2l-b)(2l-2b)(2l-3b)}{24b^{4}}t^{4}\Bigg{]}\\&\geq
b-\frac{l(2l+b)}{b}\Bigg{(}\frac{b+2(l-1)}{12(b+2l)}(bA)^{-\frac{2}{b}}\Bigg{)}^{2}
\\&~\quad-\frac{(2l-b)(8l^{2}+4lb)}{6b^{2}}\Bigg{(}\frac{b+2(l-1)}{12(b+2l)}(bA)^{-\frac{2}{b}}\Bigg{)}^{3}
\\&~\quad-\frac{(2l-b)(2l-2b)(12l^{2}+6lb)}{24b^{3}}\Bigg{(}\frac{b+2(l-1)}{12(b+2l)}(bA)^{-\frac{2}{b}}\Bigg{)}^{4}
\end{aligned}\end{eqnarray*}
 and
\begin{eqnarray*}\begin{aligned}&~\quad\Bigg{(}1+\frac{b+2(l-1)}{12(b+2l)}(bA)^{-\frac{2}{b}}\Bigg{)}^{\frac{b+2(l-1)}{b}}
\\&=(1+t)^{\frac{b+2(l-1)}{b}}\\&\geq1+\frac{2(l-1)+b}{b}\Bigg{(}\frac{b+2(l-1)}{12(b+2l)}(bA)^{-\frac{2}{b}}\Bigg{)}
\\&~\quad+\frac{(2(l-1)+b)(l-1)}{b^{2}}\Bigg{(}\frac{b+2(l-1)}{12(b+2l)}(bA)^{-\frac{2}{b}}\Bigg{)}^{2}
\\&~\quad+\frac{(2(l-1)-b)(l-1)(2(l-1)+b)}{3b^{3}}\Bigg{(}\frac{b+2(l-1)}{12(b+2l)}(bA)^{-\frac{2}{b}}\Bigg{)}^{3}.\end{aligned}\end{eqnarray*}

Therefore, we have
\begin{eqnarray*}\begin{aligned}f(\tau)&=(b+2l)\tau^{2l}bA-2l\tau^{b+2l}+\frac{l}{6}\tau^{b+2(l-1)}
\\&\geq
(bA)^{\frac{b+2l}{b}}\Bigg{[}b-\frac{l(2l+b)}{b}\Bigg{(}\frac{b+2(l-1)}{12(b+2l)}(bA)^{-\frac{2}{b}}\Bigg{)}^{2}\\&
~\quad-\frac{(2l-b)(8l^{2}+4lb)}{6b^{2}}\Bigg{(}\frac{b+2(l-1)}{12(b+2l)}(bA)^{-\frac{2}{b}}\Bigg{)}^{3}
\\&~\quad-\frac{(2l-b)(2l-2b)(12l^{2}+6lb)}{24b^{3}}\Bigg{(}\frac{b+2(l-1)}{12(b+2l)}(bA)^{-\frac{2}{b}}\Bigg{)}^{4}\Bigg{]}\\&~\quad+
\frac{l}{6}(bA)^{\frac{b+2(l-1)}{b}}\Bigg{[}1+\frac{2(l-1)+b}{b}\Bigg{(}\frac{b+2(l-1)}{12(b+2l)}(bA)^{-\frac{2}{b}}\Bigg{)}
\\&~\quad+\frac{(2(l-1)+b)(l-1)}{b^{2}}\Bigg{(}\frac{b+2(l-1)}{12(b+2l)}(bA)^{-\frac{2}{b}}\Bigg{)}^{2}
\\&~\quad+\frac{(2(l-1)-b)(l-1)(2(l-1)+b)}{3b^{3}}\Bigg{(}\frac{b+2(l-1)}{12(b+2l)}(bA)^{-\frac{2}{b}}\Bigg{)}^{3}\Bigg{]}
\\&=b(bA)^{\frac{b+2l}{b}}+\frac{l}{6}(bA)^{\frac{b+2(l-1)}{b}}+\frac{l(b+2(l-1))^{2}}{144b(b+2l)}(bA)^{\frac{b+2l-4}{b}}
+\eta_{1},\end{aligned}\end{eqnarray*}
where
\begin{eqnarray*}\begin{aligned}\eta_{1}&=\frac{2l(l+b-3)(b+2l)}{3b^{2}}\Bigg{(}\frac{b+2(l-1)}{12(b+2l)}\Bigg{)}^{3}(bA)^{\frac{b+2l-6}{b}}
\\&~\quad+\frac{l(b+2(l-1))(4(l-1)(2(l-1)-b)-3(2l-b)(l-b))}{72b^{3}}
\\&~\quad\times\Bigg{(}\frac{b+2(l-1)}{12(b+2l)}\Bigg{)}^{3}(bA)^{\frac{b+2l-8}{b}}.\end{aligned}\end{eqnarray*}
Since
$(bA)^{\frac{2}{b}}\geq\frac{1}{(b+1)^{\frac{2}{b}}}\geq\frac{1}{3}>\frac{1}{4}$
and $b\geq2l$, we have
\begin{eqnarray*}\begin{aligned}\eta_{1}&\geq\frac{2l(l+b-3)(b+2l)}{12b^{2}}\Bigg{(}\frac{b+2(l-1)}{12(b+2l)}\Bigg{)}^{3}(bA)^{\frac{b+2l-8}{b}}
\\&~\quad+\frac{l(b+2(l-1))(4(l-1)(2(l-1)-b)-3(2l-b)(l-b))}{72b^{3}}
\\&~\quad\times\Bigg{(}\frac{b+2(l-1)}{12(b+2l)}\Bigg{)}^{3}(bA)^{\frac{b+2l-8}{b}}
\\&=\frac{l\big{[}9b^{3}+(35l-26)b^{2}+(36l^{2}-90l)b+(4l^{3}-36l^{2}+48l-16)\big{]}}{72b^{3}}
\\&~\quad\times\Bigg{(}\frac{b+2(l-1)}{12(b+2l)}\Bigg{)}^{3}(bA)^{\frac{b+2l-8}{b}}
\\&\geq\frac{l\big{[}72l^{3}+(70l^{2}-52l)b+(36l^{2}-90l)b-36l^{2}\big{]}}{72b^{3}}
\Bigg{(}\frac{b+2(l-1)}{12(b+2l)}\Bigg{)}^{3}(bA)^{\frac{b+2l-8}{b}}
\\&\geq\frac{l\big{[}(72l^{3}-36l^{2})+(140l-52l)b+(72l-90l)b\big{]}}{72b^{3}}
\Bigg{(}\frac{b+2(l-1)}{12(b+2l)}\Bigg{)}^{3}(bA)^{\frac{b+2l-8}{b}}
\\&\geq0,\end{aligned}\end{eqnarray*}

which implies $$f(\tau)\geq
b(bA)^{\frac{b+2l}{b}}+\frac{l}{6}(bA)^{\frac{b+2(l-1)}{b}}+\frac{l(b+2(l-1))^{2}}{144b(b+2l)}(bA)^{\frac{b+2l-4}{b}}
.$$

\noindent\textbf{Case 2:} $\textbf{2\emph{l}}-\textbf{2}\leq
\textbf{\emph{b}}<\textbf{2\emph{l}.}$ By using Taylor formula, we
obtain the following inequalities for $t>0$:
\begin{eqnarray*}(1+t)^{\frac{2l}{b}}\geq1+\frac{2l}{b}t+\frac{2l(2l-b)}{2b^{2}}t^{2}+
\frac{2l(2l-b)(2l-2b)}{6b^{3}}t^{3}
\end{eqnarray*}
and
\begin{eqnarray*}\begin{aligned}(1+t)^{\frac{b+2(l-1)}{b}}&\geq1+\frac{2(l-1)+b}{b}t+\frac{(2(l-1)+b)(l-1)}{b^{2}}t^{2}
\\&~\quad+\frac{(2(l-1)+b)(l-1)(2(l-1)-b)}{3b^{3}}t^{3}.\end{aligned}\end{eqnarray*}
Putting
$$t=\frac{b+2(l-1)}{12(b+2l)}(bA)^{-\frac{2}{b}},$$
we have $b-2lt>\frac{l}{3}>0$,
\begin{eqnarray*}\begin{aligned}&\quad~\Bigg{(}b-\frac{l(b+2(l-1))}{6(b+2l)}(bA)^{-\frac{2}{b}}\Bigg{)}
\Bigg{(}1+\frac{b+2(l-1)}{12(b+2l)}(bA)^{-\frac{2}{b}}\Bigg{)}^{\frac{2l}{b}}\\&=
(b-2lt)(1+t)^{\frac{2l}{b}}\\&\geq(b-2lt)\Bigg{[}1+\frac{2l}{b}t+\frac{2l(2l-b)}{2b^{2}}t^{2}+
\frac{2l(2l-b)(2l-2b)}{6b^{3}}t^{3}\Bigg{]}\\&=
b-\frac{l(b+2l)}{b}\Bigg{(}\frac{b+2(l-1)}{12(b+2l)}(bA)^{-\frac{2}{b}}\Bigg{)}^{2}-\frac{(2l-b)(8l^{2}+4lb)}{6b^{2}}
\Bigg{(}\frac{b+2(l-1)}{12(b+2l)}(bA)^{-\frac{2}{b}}\Bigg{)}^{3}\\&\quad~-\frac{4l^{2}(2l-b)(2l-2b)}{6b^{3}}
\Bigg{(}\frac{b+2(l-1)}{12(b+2l)}(bA)^{-\frac{2}{b}}\Bigg{)}^{4}
\end{aligned}\end{eqnarray*}
and
\begin{eqnarray*}\begin{aligned}&~\quad\Bigg{(}1+\frac{b+2(l-1)}{12(b+2l)}(bA)^{-\frac{2}{b}}\Bigg{)}^{\frac{b+2(l-1)}{b}}
\\&=(1+t)^{\frac{b+2(l-1)}{b}}\\&\geq1+\frac{2(l-1)+b}{b}\Bigg{(}\frac{b+2(l-1)}{12(b+2l)}(bA)^{-\frac{2}{b}}\Bigg{)}
\\&~\quad+\frac{(2(l-1)+b)(l-1)}{b^{2}}\Bigg{(}\frac{b+2(l-1)}{12(b+2l)}(bA)^{-\frac{2}{b}}\Bigg{)}^{2}
\\&~\quad+\frac{(2(l-1)+b)(l-1)(2(l-1)-b)}{3b^{3}}\Bigg{(}\frac{b+2(l-1)}{12(b+2l)}(bA)^{-\frac{2}{b}}\Bigg{)}^{3}.\end{aligned}\end{eqnarray*}
Furthermore, we deduce by using the same method as the Case (1)
\begin{eqnarray*}\begin{aligned}f(\tau)&=(b+2l)\tau^{2l}bA-2l\tau^{b+2l}+\frac{l}{6}\tau^{b+2(l-1)}
\\&\geq
(bA)^{\frac{b+2l}{b}}\Bigg{[}b-\frac{l(b+2l)}{b}\Bigg{(}\frac{b+2(l-1)}{12(b+2l)}(bA)^{-\frac{2}{b}}\Bigg{)}^{2}
\\&~\quad-\frac{(2l-b)(8l^{2}+4lb)}{6b^{2}}
\Bigg{(}\frac{b+2(l-1)}{12(b+2l)}(bA)^{-\frac{2}{b}}\Bigg{)}^{3}
\\&~\quad-\frac{4l^{2}(2l-b)(2l-2b)}{6b^{3}}\Bigg{(}\frac{b+2(l-1)}{12(b+2l)}(bA)^{-\frac{2}{b}}\Bigg{)}^{4}\Bigg{]}
\\&~\quad+
\frac{l}{6}(bA)^{\frac{b+2(l-1)}{b}}\Bigg{[}1+\frac{2(l-1)+b}{b}\Bigg{(}\frac{b+2(l-1)}{12(b+2l)}(bA)^{-\frac{2}{b}}\Bigg{)}
\\&~\quad+\frac{(2(l-1)+b)(l-1)}{b^{2}}\Bigg{(}\frac{b+2(l-1)}{12(b+2l)}(bA)^{-\frac{2}{b}}\Bigg{)}^{2}
\\&~\quad+\frac{(2(l-1)+b)(l-1)(2(l-1)-b)}{3b^{3}}\Bigg{(}\frac{b+2(l-1)}{12(b+2l)}(bA)^{-\frac{2}{b}}\Bigg{)}^{3}\Bigg{]}
\\&=b(bA)^{\frac{b+2l}{b}}+\frac{l}{6}(bA)^{\frac{b+2(l-1)}{b}}+\frac{l(b+2(l-1))^{2}}{144b(b+2l)}(bA)^{\frac{b+2l-4}{b}}
+\eta_{2},\end{aligned}\end{eqnarray*} where
\begin{eqnarray*}\begin{aligned}\eta_{2}&=\frac{2l(l+b-3)(b+2l)}{3b^{2}}\Bigg{(}\frac{b+2(l-1)}{12(b+2l)}\Bigg{)}^{3}(bA)^{\frac{b+2l-6}{b}}
\\&~\quad+\frac{l(2(l-1)+b)[(l-1)(2(l-1)-b)(b+2l)-l(2l-b)(2l-2b)]}{18b^{3}(b+2l)}
\\&~\quad\times\Bigg{(}\frac{b+2(l-1)}{12(b+2l)}\Bigg{)}^{3}(bA)^{\frac{b+2l-8}{b}}
\\&\geq\frac{2l(l+b-3)(b+2l)}{9b^{2}}\Bigg{(}\frac{b+2(l-1)}{12(b+2l)}\Bigg{)}^{3}(bA)^{\frac{b+2l-8}{b}}
\\&~\quad+\frac{l(2(l-1)+b)(l-1)(2(l-1)-b)}{18b^{3}}
\Bigg{(}\frac{b+2(l-1)}{12(b+2l)}\Bigg{)}^{3}(bA)^{\frac{b+2l-8}{b}}
\\&=\Bigg{[}\frac{4bl(l+b-3)(b+2l)}{18b^{3}}
+\frac{l(2(l-1)+b)(l-1)(2(l-1)-b)}{18b^{3}}
\Bigg{]}
\\&~\quad\times\Bigg{(}\frac{b+2(l-1)}{12(b+2l)}\Bigg{)}^{3}(bA)^{\frac{b+2l-8}{b}}
\\&\geq\Bigg{[}\frac{4bl(l+b-3)(b+2l)}{18b^{3}}+\frac{lb(b+2l)(2(l-1)-b)}{18b^{3}} \Bigg{]}
\\&~\quad\times\Bigg{(}\frac{b+2(l-1)}{12(b+2l)}\Bigg{)}^{3}(bA)^{\frac{b+2l-8}{b}}
\\&\geq\frac{bl(b+2l)(6l+3b-14)}{18b^{3}}\Bigg{(}\frac{b+2(l-1)}{12(b+2l)}\Bigg{)}^{3}(bA)^{\frac{b+2l-8}{b}}
\\&\geq0
\end{aligned}\end{eqnarray*}
 since
$(bA)^{\frac{2}{b}}\geq\frac{1}{(b+1)^{\frac{2}{b}}}\geq\frac{1}{3}$.
Therefore, we have $$f(\tau)\geq
b(bA)^{\frac{b+2l}{b}}+\frac{l}{6}(bA)^{\frac{b+2(l-1)}{b}}+\frac{l(b+2(l-1))^{2}}{144b(b+2l)}(bA)^{\frac{b+2l-4}{b}}
.$$

\noindent\textbf{Case 3:} $\textbf{\emph{l}}\leq
\textbf{\emph{b}}<\textbf{2\emph{l}}-\textbf{2.}$ By using the
Taylor formula, one has for $t>0$
\begin{eqnarray*}(1+t)^{\frac{2l}{b}}\geq1+\frac{2l}{b}t+\frac{l(2l-b)}{b^{2}}t^{2}+
\frac{l(2l-b)(2l-2b)}{3b^{3}}t^{3}
\end{eqnarray*}
and
\begin{eqnarray*}(1+t)^{\frac{b+2(l-1)}{b}}\geq
1+\frac{2(l-1)+b}{b}t+\frac{(2(l-1)+b)(l-1)}{b^{2}}t^{2}.\end{eqnarray*}
Putting
$$t=\frac{b+2(l-1)}{12(b+2l)}(bA)^{-\frac{2}{b}}>0,$$ one has $b-2lt>0$,

\begin{eqnarray*}\begin{aligned}&\quad~\Bigg{(}b-\frac{l(b+2(l-1))}{6(b+2l)}(bA)^{-\frac{2}{b}}\Bigg{)}
\Bigg{(}1+\frac{b+2(l-1)}{12(b+2l)}(bA)^{-\frac{2}{b}}\Bigg{)}^{\frac{2l}{b}}
\\&=
(b-2lt)(1+t)^{\frac{2l}{b}}\\&\geq(b-2lt)\Bigg{[}1+\frac{2l}{b}t+\frac{l(2l-b)}{b^{2}}t^{2}+
\frac{l(2l-b)(2l-2b)}{3b^{3}}t^{3}\Bigg{]}\\&=
b-\frac{l(b+2l)}{b}\Bigg{(}\frac{b+2(l-1)}{12(b+2l)}(bA)^{-\frac{2}{b}}\Bigg{)}^{2}-\frac{(2l-b)(4l^{2}+2lb)}{3b^{2}}
\Bigg{(}\frac{b+2(l-1)}{12(b+2l)}(bA)^{-\frac{2}{b}}\Bigg{)}^{3}
\\&\quad~-\frac{4l^{2}(2l-b)(l-b)}{3b^{3}}\Bigg{(}\frac{b+2(l-1)}{12(b+2l)}(bA)^{-\frac{2}{b}}\Bigg{)}^{4}
\end{aligned}\end{eqnarray*}
and
\begin{eqnarray*}\begin{aligned}\Bigg{(}1+\frac{b+2(l-1)}{12(b+2l)}(bA)^{-\frac{2}{b}}\Bigg{)}^{\frac{b+2(l-1)}{b}}
&=(1+t)^{\frac{b+2(l-1)}{b}}\\&\geq
1+\frac{2(l-1)+b}{b}\frac{b+2(l-1)}{12(b+2l)}(bA)^{-\frac{2}{b}}\\&\quad~+\frac{(2(l-1)+b)(l-1)}{b^{2}}\Bigg{(}\frac{b+2(l-1)}{12(b+2l)}(bA)^{-\frac{2}{b}}\Bigg{)}^{2}
.\end{aligned}
\end{eqnarray*}

By the same argument as the Case 2, we can deduce the following
\begin{eqnarray*}\begin{aligned}f(\tau)&=(b+2l)\tau^{2l}bA-2l\tau^{b+2l}+\frac{l}{6}\tau^{b+2(l-1)}
\\&\geq
(bA)^{\frac{b+2l}{b}}\Bigg{[}b-\frac{l(b+2l)}{b}\Bigg{(}\frac{b+2(l-1)}{12(b+2l)}(bA)^{-\frac{2}{b}}\Bigg{)}^{2}
\\&\quad~-\frac{(2l-b)(4l^{2}+2lb)}{3b^{2}}
\Bigg{(}\frac{b+2(l-1)}{12(b+2l)}(bA)^{-\frac{2}{b}}\Bigg{)}^{3}
\\&\quad~-\frac{4l^{2}(2l-b)(l-b)}{3b^{3}}\Bigg{(}\frac{b+2(l-1)}{12(b+2l)}(bA)^{-\frac{2}{b}}\Bigg{)}^{4}\Bigg{]}
\\&\quad~+
\frac{l}{6}(bA)^{\frac{b+2(l-1)}{b}}\Bigg{[}1+\frac{2(l-1)+b}{b}\frac{b+2(l-1)}{12(b+2l)}(bA)^{-\frac{2}{b}}
\\&\quad~+\frac{(2(l-1)+b)(l-1)}{b^{2}}\Bigg{(}\frac{b+2(l-1)}{12(b+2l)}(bA)^{-\frac{2}{b}}\Bigg{)}^{2}\Bigg{]}
\\&=b(bA)^{\frac{b+2l}{b}}+\frac{l}{6}(bA)^{\frac{b+2(l-1)}{b}}+\frac{l(b+2(l-1))^{2}}{144b(b+2l)}(bA)^{\frac{b+2l-4}{b}}
+\eta_{3},
\end{aligned}\end{eqnarray*}where
\begin{eqnarray*}\begin{aligned}\eta_{3}&=\frac{2l(l+b-3)(b+2l)}{3b^{2}}\Bigg{(}\frac{b+2(l-1)}{12(b+2l)}\Bigg{)}^{3}(bA)^{\frac{b+2l-6}{b}}
\\&~\quad-\frac{4l^{2}(2l-b)(l-b)}{3b^{3}}\Bigg{(}\frac{b+2(l-1)}{12(b+2l)}\Bigg{)}^{4}(bA)^{\frac{b+2l-8}{b}}
\\&\geq0.\end{aligned}\end{eqnarray*} Therefore, we have
\begin{eqnarray*}f(\tau)\geq b(bA)^{\frac{b+2l}{b}}+\frac{l}{6}(bA)^{\frac{b+2(l-1)}{b}}+\frac{l(b+2(l-1))^{2}}{144b(b+2l)}(bA)^{\frac{b+2l-4}{b}}.\end{eqnarray*}

\noindent\textbf{Case 4:} $\textbf{2}\leq
\textbf{\emph{b}}<\textbf{\emph{l}.}$ Since $2\leq b<l$, there
exists a positive integer $k$ such that $2\leq
k-1\leq\frac{2l}{b}<k$, then we have for $t>0$ that
\begin{eqnarray*}\begin{aligned}(1+t)^{\frac{2l}{b}}&\geq1+\frac{2l}{b}t+\frac{1}{2!}\frac{2l}{b}
\Bigg{(}\frac{2l}{b}-1\Bigg{)}t^{2}+\frac{1}{3!}\frac{2l}{b}\Bigg{(}\frac{2l}{b}-1\Bigg{)}\Bigg{(}\frac{2l}{b}-2\Bigg{)}t^{3}
\\&~\quad+\cdots+\frac{1}{(k+1)!}\frac{2l}{b}\Bigg{(}\frac{2l}{b}-1\Bigg{)}\cdots\Bigg{(}\frac{2l}{b}-k\Bigg{)}t^{k+1}
\\&=1+\sum_{p=0}^{k}\Bigg{\{}\frac{1}{(p+1)!}\prod^{p}_{q=0}\Bigg{(}\frac{2l}{b}-q\Bigg{)}\Bigg{\}}t^{p+1},\end{aligned}\end{eqnarray*}
\begin{eqnarray*}\begin{aligned}(1+t)^{\frac{b+2l}{b}}&\leq1+\frac{b+2l}{b}t+\frac{1}{2!}\frac{b+2l}{b}\frac{2l}{b}t^{2}
+\frac{1}{3!}\frac{b+2l}{b}\frac{2l}{b}\Bigg{(}\frac{2l}{b}-1\Bigg{)}t^{3}
\\&~\quad+\cdots+\frac{1}{(k+1)!}\frac{b+2l}{b}\frac{2l}{b}\Bigg{(}\frac{2l}{b}-1\Bigg{)}\cdots\Bigg{(}\frac{2l}{b}-(k-1)\Bigg{)}t^{k+1}
\\&=1+\sum_{p=0}^{k}\Bigg{\{}\frac{1}{(p+1)!}\prod^{p}_{q=0}\Bigg{(}\frac{2l}{b}-q+1\Bigg{)}\Bigg{\}}t^{p+1},
\end{aligned}\end{eqnarray*} and
\begin{eqnarray*}\begin{aligned}(1+t)^{\frac{b+2(l-1)}{b}}&\geq1+\frac{2(l-1)+b}{b}t+\frac{1}{2!}\frac{(2(l-1)+b)}{b}\frac{2(l-1)}{b}t^{2}
\\&~\quad+\frac{1}{3!}\frac{(2(l-1)+b)}{b}\frac{2(l-1)}{b}\Bigg{(}\frac{2(l-1)}{b}-1\Bigg{)}t^{3}
\\&~\quad+\cdots+\frac{1}{k!}\frac{(2(l-1)+b)}{b}\frac{2(l-1)}{b}\cdots\Bigg{(}\frac{2(l-1)}{b}-(k-2)\Bigg{)}t^{k}
\\&~\quad-\Bigg{|}\frac{1}{(k+1)!}\frac{(2(l-1)+b)}{b}\frac{2(l-1)}{b}\cdots\Bigg{(}\frac{2(l-1)}{b}-(k-1)\Bigg{)}\Bigg{|}t^{k+1}
\\&=1+\sum_{p=0}^{k-1}\Bigg{\{}\frac{1}{(p+1)!}\prod^{p}_{q=0}\Bigg{(}\frac{2(l-1)}{b}-q+1\Bigg{)}\Bigg{\}}t^{p+1}
\\&~\quad-\Bigg{|}\frac{1}{(k+1)!}\prod^{k}_{q=0}\Bigg{(}\frac{2(l-1)}{b}-q+1\Bigg{)}\Bigg{|}t^{k+1}.\end{aligned}\end{eqnarray*}
Putting $t=\frac{b+2(l-1)}{12(b+2l)}(bA)^{-\frac{2}{b}}$ and
$f(\tau)=(bA)^{\frac{b+2l}{b}}h(\tau),$ where
$$h(\tau)=(b+2l)(1+t)^{\frac{2l}{b}}-2l(1+t)^{\frac{b+2l}{b}}+\frac{1}{6}(bA)^{-\frac{2}{b}}(1+t)^{\frac{b+2(l-1)}{b}},$$
then we have for $2\leq b<l$,
\begin{eqnarray*}\begin{aligned}h(\tau)&\geq(b+2l)\Bigg{\{}1+\sum_{p=0}^{k}
\Bigg{[}\frac{1}{(p+1)!}\prod^{p}_{q=0}\Bigg{(}\frac{2l}{b}-q\Bigg{)}\Bigg{]}t^{p+1}\Bigg{\}}
\\&~\quad-2l\Bigg{\{}1+\sum_{p=0}^{k}\Bigg{[}\frac{1}{(p+1)!}\prod^{p}_{q=0}\Bigg{(}\frac{2l}{b}-q+1\Bigg{)}\Bigg{]}t^{p+1}\Bigg{\}}
\\&~\quad+\frac{l}{6}(bA)^{-\frac{2}{b}}\Bigg{\{}1+\sum_{p=0}^{k-1}\Bigg{[}\frac{1}{(p+1)!}\prod^{p}_{q=0}\Bigg{(}\frac{2(l-1)}{b}-q+1\Bigg{)}\Bigg{]}t^{p+1}
\\&~\quad-\Bigg{|}\frac{1}{(k+1)!}\prod^{k}_{q=0}\Bigg{(}\frac{2(l-1)}{b}-q+1\Bigg{)}\Bigg{|}t^{k+1}\Bigg{\}}
\\&=b+\frac{l}{6}(bA)^{-\frac{2}{b}}+\sum_{p=1}^{k}\Bigg{\{}\frac{b+2l}{(p+1)!}\frac{2l}{b}\Bigg{[}\prod^{p}_{q=1}
\Bigg{(}\frac{2l}{b}-q\Bigg{)}
-\prod^{p}_{q=1}\Bigg{(}\frac{2l}{b}-q+1\Bigg{)}\Bigg{]}\Bigg{\}}t^{p+1}
\\&~\quad+\sum_{p=0}^{k-1}\Bigg{\{}\frac{l(bA)^{-\frac{2}{b}}}{6(p+1)!}\prod^{p}_{q=0}
\Bigg{(}\frac{2(l-1)}{b}-q+1\Bigg{)}\Bigg{\}}t^{p+1}
\\&~\quad-\Bigg{|}\frac{l(bA)^{-\frac{2}{b}}}{6(k+1)!}\prod^{k}_{q=0}
\Bigg{(}\frac{2(l-1)}{b}-q+1\Bigg{)}\Bigg{|}t^{k+1}
\\&=b+\frac{l}{6}(bA)^{-\frac{2}{b}}-\sum_{p=1}^{k}\Bigg{\{}\frac{p2l(b+2l)}{b(p+1)!}\prod^{p-1}_{q=1}
\Bigg{(}\frac{2l}{b}-q\Bigg{)}\Bigg{\}}t^{p+1}
\\&~\quad+\sum_{p=1}^{k}\Bigg{\{}\frac{l(bA)^{-\frac{2}{b}}}{6p!}\prod^{p-1}_{q=0}
\Bigg{(}\frac{2(l-1)}{b}-q+1\Bigg{)}\Bigg{\}}t^{p}
\\&~\quad-\Bigg{|}\frac{l(bA)^{-\frac{2}{b}}}{6(k+1)!}\prod^{k}_{q=0}
\Bigg{(}\frac{2(l-1)}{b}-q+1\Bigg{)}\Bigg{|}t^{k+1}
\\&=
b+\frac{l}{6}(bA)^{-\frac{2}{b}}-\sum_{p=1}^{k}\Bigg{\{}\frac{p}{b^{p}(p+1)!}\prod^{p}_{q=0}(2l-(q-1)b)\Bigg{\}}
\Bigg{(}\frac{b+2(l-1)}{12(b+2l)}(bA)^{-\frac{2}{b}}\Bigg{)}^{p+1}
\\&~\quad+\sum_{p=1}^{k}\Bigg{\{}\frac{l(bA)^{-\frac{2}{b}}}{6b^{p}p!}\prod^{p-1}_{q=0}(2(l-1)-(q-1)b)\Bigg{\}}
\Bigg{(}\frac{b+2(l-1)}{12(b+2l)}(bA)^{-\frac{2}{b}}\Bigg{)}^{p}
\\&~\quad-\Bigg{|}\frac{l(bA)^{-\frac{2}{b}}}{6b^{k+1}(k+1)!}\prod^{k}_{q=0}(2(l-1)-(q-1)b)
\Bigg{|}\Bigg{(}\frac{b+2(l-1)}{12(b+2l)}(bA)^{-\frac{2}{b}}\Bigg{)}^{k+1}.\end{aligned}\end{eqnarray*}

Furthermore,
\begin{eqnarray*}\begin{aligned}f(\tau)&\geq
b(bA)^{\frac{b+2l}{b}}+\frac{l}{6}(bA)^{\frac{b+2(l-1)}{b}}+\frac{l(b+2(l-1))^{2}}{144b(b+2l)}(bA)^{\frac{b+2l-4}{b}}
\\&~\quad-\sum_{p=2}^{k}\Bigg{\{}\frac{p}{b^{p}(p+1)!}\prod^{p}_{q=0}(2l-(q-1)b)\Bigg{\}}
\Bigg{(}\frac{b+2(l-1)}{12(b+2l)}\Bigg{)}^{p+1}(bA)^{\frac{b+2l-2p-2}{b}}\end{aligned}\end{eqnarray*}

\begin{eqnarray*}\begin{aligned}
\quad\quad~&~\quad+\sum_{p=2}^{k}\Bigg{\{}\frac{l}{6b^{p}p!}\prod^{p-1}_{q=0}(2(l-1)-(q-1)b)\Bigg{\}}
\Bigg{(}\frac{b+2(l-1)}{12(b+2l)}\Bigg{)}^{p}(bA)^{\frac{b+2l-2p-2}{b}}
\\&~\quad-\Bigg{|}\frac{l}{6b^{k+1}(k+1)!}\prod^{k}_{q=0}(2(l-1)-(q-1)b)
\Bigg{|}\Bigg{(}\frac{b+2(l-1)}{12(b+2l)}\Bigg{)}^{k+1}(bA)^{\frac{b+2l-2k-4}{b}}
\\&=
b(bA)^{\frac{b+2l}{b}}+\frac{l}{6}(bA)^{\frac{b+2(l-1)}{b}}+\frac{l(b+2(l-1))^{2}}{144b(b+2l)}(bA)^{\frac{b+2l-4}{b}}
+\eta_{4} ,\end{aligned}\end{eqnarray*}
where
\begin{eqnarray*}\begin{aligned}\eta_{4}&=\sum_{p=2}^{k}
\Bigg{\{}\frac{(b+2(l-1))2l}{12b^{p}p!}\Bigg{[}\prod^{p-1}_{q=1}(2(l-1)-(q-1)b)
-\frac{p}{p+1}\prod^{p-1}_{q=1}(2l-qb)\Bigg{]}\Bigg{\}}
\\&~\quad\times
\Bigg{(}\frac{b+2(l-1)}{12(b+2l)}\Bigg{)}^{p}(bA)^{\frac{b+2l-2p-2}{b}}\\&~\quad-\Bigg{|}\frac{l}{6b^{k+1}(k+1)!}\prod^{k}_{q=0}(2(l-1)-(q-1)b)
\Bigg{|}\Bigg{(}\frac{b+2(l-1)}{12(b+2l)}\Bigg{)}^{k+1}(bA)^{\frac{b+2l-2k-4}{b}}.
\end{aligned}\end{eqnarray*}
From $k-2\leq\frac{2(l-1)}{b}<k,$ we have
\begin{eqnarray}\label{2.7}\frac{k-2-i}{k+1-i}\leq\frac{\frac{2(l-1)}{b}-i}{k+1-i}<\frac{k-i}{k+1-i}.\end{eqnarray}
Then, it follows from (\ref{2.7}) that
\begin{eqnarray*}
\Bigg{|}\frac{\frac{2(l-1)}{b}-i}{k+1-i}\Bigg{|}\leq1, \ \ \ {\rm for}\ i=0,1,2,\cdots,k-1.
\end{eqnarray*}
Note that
$$2(l-1)-(q-1)b\geq2l-qb\geq0, \ \ {\rm for}\ p=2,3,\cdots,k,$$
one has
\begin{eqnarray*}
\begin{aligned}
&\quad~\prod^{p-1}_{q=1}(2(l-1)-(q-1)b)
-\frac{p}{p+1}\prod^{p-1}_{q=1}(2l-qb)\\
&\geq\prod^{p-1}_{q=1}(2(l-1)-(q-1)b)
-\prod^{p-1}_{q=1}(2l-qb)\geq 0.
\end{aligned}
\end{eqnarray*}
Therefore, we obtain
\begin{eqnarray*}
\begin{aligned}&\sum_{p=2}^{k}\Bigg{\{}\frac{(b+2(l-1))2l}{12b^{p}p!}\Bigg{[}\prod^{p-1}_{q=1}(2(l-1)-(q-1)b)
-\frac{p}{p+1}\prod^{p-1}_{q=1}(2l-qb)\Bigg{]}\Bigg{\}}
\\&~\quad\times
\Bigg{(}\frac{b+2(l-1)}{12(b+2l)}\Bigg{)}^{p}(bA)^{\frac{b+2l-2p-2}{b}}\\&\geq
\frac{(b+2(l-1))2l}{24b^{2}}\Bigg{[}2(l-1)-\frac{2(2l-b)}{3}\Bigg{]}
\Bigg{(}\frac{b+2(l-1)}{12(b+2l)}\Bigg{)}^{2}(bA)^{\frac{b+2l-6}{b}}.
\end{aligned}
\end{eqnarray*}
From
$$(bA)^{\frac{2}{b}}\geq\frac{1}{(b+1)^{\frac{2}{b}}}\geq\frac{1}{3},$$
we have
\begin{eqnarray*}\begin{aligned}\eta_{4}&\geq\frac{(b+2(l-1))2l}{24b^{2}}\Bigg{[}2(l-1)-\frac{2(2l-b)}{3}\Bigg{]}
\Bigg{(}\frac{b+2(l-1)}{12(b+2l)}\Bigg{)}^{2}(bA)^{\frac{b+2l-6}{b}}
\\&~\quad-\Bigg{|}\frac{l}{6b^{k+1}(k+1)!}\prod^{k}_{q=0}(2(l-1)-(q-1)b)
\Bigg{|}\Bigg{(}\frac{b+2(l-1)}{12(b+2l)}\Bigg{)}^{k+1}(bA)^{\frac{b+2l-2k-4}{b}}
\\&=\frac{l(b+2(l-1))}{12b^{2}}
\Bigg{\{}\Bigg{[}2(l-1)-\frac{2(2l-b)}{3}\Bigg{]}
\\&~\quad-2b\Bigg{|}\frac{1}{b^{k}(k+1)!}\prod^{k}_{q=1}(2(l-1)-(q-1)b)
\Bigg{|}\Bigg{(}\frac{b+2(l-1)}{12(b+2l)}\Bigg{)}^{k-1}(bA)^{\frac{-2k+2}{b}}\Bigg{\}}
\\&~\quad\times\Bigg{(}\frac{b+2(l-1)}{12(b+2l)}\Bigg{)}^{2}(bA)^{\frac{b+2l-6}{b}}\\&\geq\frac{l(b+2(l-1))}{12b^{2}}
\Bigg{\{}\frac{2l+2b-6}{3}
-2b\Bigg{|}\prod^{k-1}_{q=0}\Bigg{(}\frac{\dfrac{2(l-1)}{b}-q}{k+1-q}\Bigg{)}
\Bigg{|}\Bigg{(}\frac{1}{4}\Bigg{)}^{k-1}\Bigg{\}}
\\&~\quad\times\Bigg{(}\frac{b+2(l-1)}{12(b+2l)}\Bigg{)}^{2}(bA)^{\frac{b+2l-6}{b}}
\\&\geq
\frac{l(b+2(l-1))}{12b^{2}}\Bigg{\{}\frac{2l+2b-6}{3}-\frac{b}{8}\Bigg{\}}\Bigg{(}\frac{b+2(l-1)}{12(b+2l)}\Bigg{)}^{2}(bA)^{\frac{b+2l-6}{b}}
\\&\geq0,\end{aligned}\end{eqnarray*}
which implies that
\begin{eqnarray*}f(\tau)\geq b(bA)^{\frac{b+2l}{b}}+\frac{l}{6}(bA)^{\frac{b+2(l-1)}{b}}+\frac{l(b+2(l-1))^{2}}{144b(b+2l)}(bA)^{\frac{b+2l-4}{b}}.\end{eqnarray*}
This completes the proof of Lemma \ref{lem2.1}.
$$\eqno{\Box}$$

\section{Proof of Theorem \ref{thm1.1} and Remark \ref{rem1.2}}\label{sec3}

\emph{Proof of Theorem \ref{thm1.1}.} We will use the same notations
as those of \cite{CQW}.   In this section, we assume that $b=n$. Let
$\widehat{\varphi}_{j}(z)$ be the Fourier transform of the trial
function $\varphi_{j}(x)$,
\begin{equation*}\varphi_{j}(x)=\begin{cases}u_{j}(x),~~~~~~~~~x\in\Omega,\\
0,~~~~~~~~x\in\mathbb{R}^{n}\setminus\Omega,\end{cases}\end{equation*}
where $u_{j}(x)$ is an orthonormal eigenfunction corresponding to
the eigenvalue $\lambda_{j}$,
$f(z):=\sum^{k}_{j=1}|\widehat{\varphi}_{j}(z)|^{2}$, and $f^{\ast}$
be the symmetric decreasing rearrangement of $f$. And then, we can
obtain from Lemma \ref{lem2.1} that
\begin{eqnarray}\label{3.1}\begin{aligned}\sum_{j=1}^{k}\lambda_{j}&\geq nB_{n}\int^{\infty}_{0}s^{n+2l-1}\phi(s)ds\\&\geq
\frac{nB_{n}\Big{(}\frac{k}{B_{n}}\Big{)}^{\frac{n+2l}{n}}}{n+2l}\phi(0)^{-\frac{2l}{n}}
+\frac{lB_{n}\Big{(}\frac{k}{B_{n}}\Big{)}^{\frac{n+2(l-1)}{n}}}{6(n+2l)\mu^{2}}\phi(0)^{\frac{2n-2l+2}{n}}
\\&~\quad+\frac{l(n+2(l-1))^{2}B_{n}\Big{(}\frac{k}{B_{n}}\Big{)}^{\frac{n+2l-4}{n}}}{144n(n+2l)^{2}\mu^{4}}
\phi(0)^{\frac{4n-2l+4}{n}},\end{aligned}\end{eqnarray} where
$\phi:[0,~+\infty)\rightarrow[0,~(2\pi)^{-n}V(\Omega)]$ is a
non-increasing function of $|x|$ and $\phi(x)$ is defined by
$\phi(|x|):=f^{\ast}(x)$. Now defining a function $\xi(t)$ as
follows:
\begin{eqnarray}\label{3.2}\begin{aligned}\xi(t)&=
\frac{nB_{n}}{n+2l}\Bigg{(}\frac{k}{B_{n}}\Bigg{)}^{\frac{n+2l}{n}}t^{-\frac{2l}{n}}
+\frac{lB_{n}}{6(n+2l)\mu^{2}}\Bigg{(}\frac{k}{B_{n}}\Bigg{)}^{\frac{n+2(l-1)}{n}}t^{\frac{2n-2l+2}{n}}
\\&~\quad+\frac{l(n+2(l-1))^{2}B_{n}}{144n(n+2l)^{2}\mu^{4}}
\Bigg{(}\frac{k}{B_{n}}\Bigg{)}^{\frac{n+2l-4}{n}}t^{\frac{4n-2l+4}{n}}.\end{aligned}\end{eqnarray}
Here we assume that $l\leq n+1$. The other cases (i.e.,
$n+1<l<2(n+1)$, $l\geq2(n+1)$) can be discussed by using of the
similar method. After differentiating (\ref{3.2}) with respect to
the variable $t$, we derive

\begin{eqnarray}\label{3.3}\begin{aligned}\xi^{\prime}(t)&=
\frac{B_{n}}{n+2l}\Bigg{(}\frac{k}{B_{n}}\Bigg{)}^{\frac{n+2l}{n}}t^{-\frac{2l}{n}-1}\Bigg{[}-2l
+\frac{l(2n-2l+2)}{6n\mu^{2}}\Bigg{(}\frac{k}{B_{n}}\Bigg{)}^{-\frac{2}{n}}t^{\frac{2n+2}{n}}
\\&~\quad+\frac{l(4n-2l+4)(n+2(l-1))^{2}}{144n^{2}(n+2l)\mu^{4}}
\Bigg{(}\frac{k}{B_{n}}\Bigg{)}^{-\frac{4}{n}}t^{\frac{4n+4}{n}}\Bigg{]}.
\end{aligned}\end{eqnarray} Putting
$\zeta(t)=\xi^{\prime}(t)\frac{n+2l}{B_{n}}(\frac{k}{B_{n}})^{-\frac{n+2l}{n}}t^{\frac{2l}{n}+1}$
and noticing that
$\mu\geq(2\pi)^{-n}B_{n}^{-\frac{1}{n}}V(\Omega)^{\frac{n+1}{n}}$,
we can deduce from (\ref{3.3})
\begin{eqnarray}\label{3.4}\begin{aligned}
\zeta(t)&= -2l
+\frac{l(2n-2l+2)}{6n\mu^{2}}\Bigg{(}\frac{k}{B_{n}}\Bigg{)}^{-\frac{2}{n}}t^{\frac{2n+2}{n}}
\\&~\quad+\frac{l(4n-2l+4)(n+2(l-1))^{2}}{144n^{2}(n+2l)\mu^{4}}
\Bigg{(}\frac{k}{B_{n}}\Bigg{)}^{-\frac{4}{n}}t^{\frac{4n+4}{n}}
\\&\leq-2l
+\frac{l(2n-2l+2)}{6n(2\pi)^{-2n}B_{n}^{-\frac{2}{n}}Vol(\Omega)^{\frac{2(n+1)}{n}}}\Bigg{(}\frac{k}{B_{n}}\Bigg{)}^{-\frac{2}{n}}t^{\frac{2n+2}{n}}
\\&~\quad+\frac{l(4n-2l+4)(n+2(l-1))^{2}}{144n^{2}(n+2l)(2\pi)^{-4n}B_{n}^{-\frac{4}{n}}Vol(\Omega)^{\frac{4(n+1)}{n}}}
\Bigg{(}\frac{k}{B_{n}}\Bigg{)}^{-\frac{4}{n}}t^{\frac{4n+4}{n}}.\end{aligned}\end{eqnarray}
Since the right hand side of (\ref{3.4}) is an increasing function
of $t$, if the right hand side of (\ref{3.4}) is not larger than $0$
at $t=(2\pi)^{-n}V(\Omega)$, that is
\begin{eqnarray}\label{3.5}\begin{aligned}
\zeta(t)&\leq-2l
+\frac{l(2n-2l+2)}{6n}k^{-\frac{2}{n}}\frac{B_{n}^{\frac{4}{n}}}{(2\pi)^{2}}
\\&~\quad+\frac{l(4n-2l+4)(n+2(l-1))^{2}}{144n^{2}(n+2l)}
k^{-\frac{4}{n}}\frac{B_{n}^{\frac{8}{n}}}{(2\pi)^{4}}
\\&\leq0,\end{aligned}\end{eqnarray}we can claim from
(\ref{3.5}) that $\xi'(t)\leq0$ on $(0,(2\pi)^{-n}V(\Omega)].$ If
$\xi'(t)\leq0$, then $\xi(t)$ is a decreasing function on
$(0,(2\pi)^{-n}V(\Omega)].$ In fact, by a direct calculation, we can
obtain

\begin{eqnarray}\label{3.6}
\zeta(t)\leq-2l +\frac{l(2n-2l+2)}{6n}
+\frac{l(4n-2l+4)(n+2(l-1))^{2}}{144n^{2}(n+2l)} \leq0\end{eqnarray}
since $\frac{B_{n}^{\frac{4}{n}}}{(2\pi)^{2}}<1.$

On the other hand, since $0<\phi(0)\leq(2\pi)^{-n}V(\Omega)$ and
right hand side of the formula (\ref{3.1}) is $\xi(\phi(0))$, which
is a decreasing function of $\phi(0)$ on $(0,(2\pi)^{-n}V(\Omega)]$,
then we can replace $\phi(0)$ by $(2\pi)^{-n}V(\Omega)$ in
(\ref{3.1}) which gives the inequality as follows:

\begin{eqnarray*}\begin{aligned}\frac{1}{k}\sum_{j=1}^{k}\lambda_{j}&\geq
\frac{n}{n+ 2l}\frac{(2\pi)^{2l}}
{(B_{n}V(\Omega))^{\frac{2l}{n}}}k^{\frac{2l}{n}}
\\&~\quad+\frac{l}{24(n+2l)}\frac{(2\pi)^{2(l-1)}}{(B_{n}V(\Omega))^{\frac{2(l-1)}{n}}}
\frac{V(\Omega)}{I(\Omega)}k^{\frac{2(l-1)}{n}}
\\&~\quad+\frac{l(n+2(l-1))^{2}}{2304n(n+2l)^{2}}
\frac{(2\pi)^{2(l-2)}}{(B_{n}V(\Omega))^{\frac{2(l-2)}{n}}}
\Bigg{(}\frac{V(\Omega)}{I(\Omega)}\Bigg{)}^{2}k^{\frac{2(l-2)}{n}}
.\end{aligned}\end{eqnarray*}

This completes the proof of Theorem \ref{thm1.1}.
$$\eqno{\Box}$$

Next we will prove that the inequality (\ref{1.17}) is sharper than
the inequality (\ref{1.16}).

\noindent\emph{Proof of Remark \ref{rem1.2}:} Under the same
assumption with Lemma \ref{lem2.1}, let $b=n$ and
$A=\frac{k}{nB_n}$, we obtain from
$\mu\geq(2\pi)^{-n}B_{n}^{-\frac{1}{n}}V(\Omega)^{\frac{n+1}{n}}$
that

\begin{eqnarray*}\begin{aligned}
\frac{(bA)^{-\frac{2}{b}}\psi(0)^{2+\frac{2}{b}}}{\mu^2}
&\leq\frac{(2\pi)^{-2n-2}V(\Omega)^{2+\frac{2}{n}}}{(2\pi)^{-2n}B_n^{-\frac{2}{n}}
V(\Omega)^{\frac{2(n+1)}{n}}}(\frac{k}{B_n})^{-\frac{2}{n}}\\
&=\frac{(2\pi)^{-2}}{(B_n)^{-\frac{4}{n}}}k^{-\frac{2}{n}}<\frac{(2\pi)^{-2}}{(B_n)^{-\frac{4}{n}}}
<1,\end{aligned}
\end{eqnarray*}
then we have
\begin{eqnarray*}\begin{aligned}
&\quad\frac{1}{b+2l}\sum_{p=2}^{l}
\frac{(l+1-p)}{(6)^{p}b\cdots(b+2p-2)\mu^{2p}} (bA)^{\frac{b+2(l-p)}{b}}\psi(0)^{\frac{2pb-2(l-p)}{b}}\\
&<\frac{1}{b+2l}\sum_{p=2}^{l}
\frac{(l+1-p)}{(6)^{p}b\cdots(b+2p-2)\mu^{4}} (bA)^{\frac{b+2l-4}{b}}\psi(0)^{\frac{4b-2l+4}{b}}\end{aligned}\end{eqnarray*}
\begin{eqnarray}\label{3.7}\begin{aligned}
&<\frac{l-1}{36 (b+2l)b(b+2)\mu^{4}}\sum_{p=0}^{\infty} \frac{1}{6^p (b+2)^p}(bA)^{\frac{b+2l-4}{b}}\psi(0)^{\frac{4b-2l+4}{b}}\\
&=\frac{l-1}{6b(b+2l)(6(b+2)-1)\mu^{4}}(bA)^{\frac{b+2l-4}{b}}\psi(0)^{\frac{4b-2l+4}{b}}.
\end{aligned}\end{eqnarray}
By a direct calculation, we derive
\begin{eqnarray*}
l(6(b+2)-1)(b+2(l-1))^{2}>24b(b+2l)(l-1)>0,
\end{eqnarray*}
in fact,
$$
\aligned &\ \ \ \
l(6(b+2)-1)(b+2(l-1))^{2}-24b(b+2l)(l-1)\\
&=4b(l-1)[6b(l-1)-l]+l(6b+11)[b^2+4(l-1)^2]\\
&>24b^2(l-1)^2+4bl(l-1)[6(l-1)-1]>0,
\endaligned
$$
that is,
\begin{eqnarray}\label{3.8}
\frac{24b(b+2l)(l-1)}{l(6(b+2)-1)(b+2(l-1))^{2}}<1.
\end{eqnarray}
Therefore, we get from (\ref{3.7}) and (\ref{3.8}) that
\begin{eqnarray}\label{3.9}\begin{aligned}
&\frac{1}{b+2l}\sum_{p=2}^{l}
\frac{(l+1-p)}{(6)^{p}b\cdots(b+2p-2)\mu^{2p}} (bA)^{\frac{b+2(l-p)}{b}}\psi(0)^{\frac{2pb-2(l-p)}{b}}\\
&<\frac{l-1}{6b(b+2l)(6(b+2)-1)\mu^{4}}(bA)^{\frac{b+2l-4}{b}}\psi(0)^{\frac{4b-2l+4}{b}}\\
&=\frac{24b(b+2l)(l-1)}{l(6(b+2)-1)(b+2(l-1))^{2}}\cdot\frac{l(b+2(l-1))^{2}}{144
b^{2}(b+2l)^{2}\mu^4}(bA)^{\frac{b+2l-4}{b}}\psi(0)^{\frac{4b-2l+4}{b}}
\\&<\frac{l(b+2(l-1))^{2}}{144
b^{2}(b+2l)^{2}\mu^4}(bA)^{\frac{b+2l-4}{b}}\psi(0)^{\frac{4b-2l+4}{b}}.
\end{aligned}\end{eqnarray}
Taking
\begin{eqnarray}\label{3.10}
b=n,\ \  A=\frac{k}{nB_n},\ \ \psi(0)=(2\pi)^{-n}V(\Omega),\ \
\mu=2(2\pi)^{-n}\sqrt{V(\Omega)I(\Omega)},
\end{eqnarray}
and substituting (\ref{3.10}) into (\ref{3.9}), one has

\begin{equation}\label{3.11}
\begin{aligned}
&
\frac{n}{n+ 2l}\frac{(2\pi)^{2l}}
{(B_{n}V(\Omega))^{\frac{2l}{n}}}k^{\frac{2l}{n}}
+\frac{l}{24(n+2l)}\frac{(2\pi)^{2(l-1)}}{(B_{n}V(\Omega))^{\frac{2(l-1)}{n}}}
\frac{V(\Omega)}{I(\Omega)}k^{\frac{2(l-1)}{n}}\\
&+\frac{l(n+2(l-1))^{2}}{2304n(n+2l)^{2}}
\frac{(2\pi)^{2(l-2)}}{(B_{n}V(\Omega))^{\frac{2(l-2)}{n}}}
\Bigg{(}\frac{V(\Omega)}{I(\Omega)}\Bigg{)}^{2}k^{\frac{2(l-2)}{n}}\\
>&\frac{n}{n+
2l}\frac{(2\pi)^{2l}}
{(B_{n}V(\Omega))^{\frac{2l}{n}}}k^{\frac{2l}{n}}+\frac{n}{(n+2l)}\\
&\times\sum_{p=1}^{l}
\frac{l+1-p}{(24)^{p}n\cdots(n+2p-2)}\frac{(2\pi)^{2(l-p)}}{(B_{n}V(\Omega))^{\frac{2(l-p)}{n}}}
\Bigg{(}\frac{V(\Omega)}{I(\Omega)}\Bigg{)}^{p}k^{\frac{2(l-p)}{n}}.
\end{aligned}
\end{equation}
This completes the proof of Remark \ref{rem1.2}.
$$\eqno{\Box}$$
\textbf{Acknowledgments}~The  authors  wish  to  express their
gratitude to Prof. Q. -M. Cheng for continuous encouragement and
enthusiastic help.

\bibliographystyle{amsplain}

\underline{}

\begin{flushleft}
\medskip\noindent

Guoxin Wei, School of Mathematical Sciences, South China Normal
University, 510631, Guangzhou, China, weigx03@mails.tsinghua.edu.cn

Lingzhong Zeng, Department of Mathematics, Graduate School of
Science and Engineering, Saga University, Saga 840-8502, Japan,
lingzhongzeng@yeah.net
\end{flushleft}

\end{document}